\pdfoutput=1
\RequirePackage{ifpdf}
\ifpdf 
\documentclass[pdftex]{sigma}
\else
\documentclass{sigma}
\fi

\numberwithin{equation}{section}

\newtheorem{Theorem}{Theorem}[section]
\newtheorem{Corollary}[Theorem]{Corollary}
\newtheorem{Lemma}[Theorem]{Lemma}
\newtheorem{Proposition}[Theorem]{Proposition}
{ \theoremstyle{definition}
\newtheorem{Definition}[Theorem]{Definition}
 }

\newcommand{\mb}{\mathbf}

\renewcommand{\O}{\mathcal{O}}

\renewcommand{\t}{\theta}

\begin{document}

\allowdisplaybreaks

\renewcommand{\thefootnote}{$\star$}

\renewcommand{\PaperNumber}{092}

\FirstPageHeading

\ShortArticleName{Orthogonal Basic Hypergeometric Laurent Polynomials}

\ArticleName{Orthogonal Basic Hypergeometric\\ Laurent Polynomials\footnote{This
paper is a contribution to the Special Issue ``Superintegrability, Exact Solvability, and Special Functions''. The full collection is available at \href{http://www.emis.de/journals/SIGMA/SESSF2012.html}{http://www.emis.de/journals/SIGMA/SESSF2012.html}}}

\Author{Mourad E.H.~ISMAIL~$^{\dag\ddag}$ and Dennis STANTON~$^\S$}

\AuthorNameForHeading{M.E.H.~Ismail and D.~Stanton}

\Address{$^\dag$~Department  of Mathematics,
  University of  Central Florida,
  Orlando, FL 32816, USA}
\EmailD{\href{mailto:mourad.eh.ismail@gmail.com}{mourad.eh.ismail@gmail.com}}

\Address{$^\ddag$~Department of Mathematics, King Saud University,
Riyadh, Saudi Arabia}

\Address{$^\S$~School of Mathematics, University of Minnesota, Minneapolis, MN 55455, USA}
\EmailD{\href{mailto:stanton@math.umn.edu}{stanton@math.umn.edu}}

\ArticleDates{Received August 04, 2012, in f\/inal form November 28, 2012; Published online December 01, 2012}

\Abstract{The Askey--Wilson polynomials are orthogonal polynomials in
$x = \cos \theta$, which
are given as a terminating $_4\phi_3$ basic hypergeometric series.
The non-symmetric Askey--Wilson polynomials are Laurent polynomials in
$z=e^{i\theta}$, which are given as a sum of two terminating $_4\phi_3$'s.
They satisfy a biorthogonality relation. In this paper new orthogo\-na\-li\-ty
relations for single  $_4\phi_3$'s which are Laurent polynomials in~$z$ are given,
which imply the non-symmetric Askey--Wilson biorthogonality. These results include
discrete ortho\-go\-na\-li\-ty relations. They can be considered as a classical analytic
study of the results for non-symmetric
Askey--Wilson polynomials which were previously obtained by af\/f\/ine Hecke
algebra techniques.}


\Classification{33D45}

\renewcommand{\thefootnote}{\arabic{footnote}}
\setcounter{footnote}{0}

\section{Introduction}

The Askey--Wilson polynomials $p_n(x;\mb{t}\,|\, q) $
are polynomials in $x = \cos \theta$, and depend
upon parameters $q$, $t_1$, $t_2$, $t_3$, and $t_4$.
They may be def\/ined by the
terminating basic hypergeometric series \cite{And:Ask:Roy}, \cite[\S~15.2]{Ism05}
\begin{gather}
p_n(x;\mb{t}\,|\, q)  =t_1^{-n}(t_1t_2, t_1t_3, t_1t_4;q)_n
\, {}_{4}\phi_3\left(\left. \begin{matrix}
q^{-n}, t_1t_2t_3t_4q^{n-1}, t_1e^{i\theta}, t_1e^{-i\theta} \\
t_1t_2,  t_1t_3,  t_1t_4
\end{matrix}  \right|q,q\right).\label{eqasexpl}
\end{gather}
They are orthogonal polynomials in~$x$, and have a known orthogonality
relation (see~\eqref{eqaworthrel}).

Clearly the Askey--Wilson polynomials  are also Laurent polynomials in $z=e^{i\theta}$, since $2x=z+1/z$.
One may ask if there is a
natural basis for Laurent polynomials in $z$ which also satisfy orthogonality relations.
The non-symmetric Askey--Wilson polynomials \cite{M,NS}
are Laurent polynomials, which are
biorthogonal with respect to a modif\/ied Askey--Wilson weight, are one answer to this question.
They may be explicitly given
as a sum of two terminating basic hypergeometric series.

The purpose of
this paper is to use classical analytic methods to study the orthogonality
relations for single ${}_{4}\phi_3$'s
which imply the non-symmetric Askey--Wilson polynomial biorthogonality.
The main orthogonality results are
Theorems~\ref{orthbasis1},~\ref{orthbasis2},~\ref{orthbasis3}, and~\ref{orthbasis4}.
We show that special function techniques, explicit summations and self-adjoint operators,
may be applied to derive these orthogonalities.
We explicitly
f\/ind three and four term relations for Laurent polynomials in Section~\ref{section4}.
As a corollary of our computations we rederive
the $L^2$-norm of the non-symmetric Askey--Wilson polynomials in Section~\ref{section6}, and
reprove the creation type
formulas for the Laurent polynomials in Section~\ref{section5}. Asymptotics are given in
Theorem~\ref{asympt}, and
Racah orthogonality in Theorem~\ref{Racah}.
This work may be considered as a classical study of the one
variable case in Cherednik~\cite{Chr}, Macdonald~\cite{M},
carefully studied by  Noumi and Stokman~\cite{NS}, who considered
many of these results from the point of view of af\/f\/ine Hecke algebras.

There are dif\/ferent theories of orthogonal polynomials, biorthogonal polynomials and rational functions
which originate from dif\/ferent types of continued fractions. The $J$-fractions lead to orthogonal
polynomials on the line while the $PC$ fractions lead to orthogonal polynomials on the unit circle~\cite{Ism05, Jon:Thr,Lor:Waa2,Lor:Waa1}. The theory of $T$-fractions leads
to orthogonal Laurent polynomials  and
is explained in the books~\cite{Lor:Waa2,Lor:Waa1},  and the excellent survey
article~\cite{Bul:Gon:Hen:Nja}. Ismail and Masson~\cite{Ism:Mas} introduced the more general $R_{\rm I}$ and $R_{\rm II}$
fractions which naturally lead to biorthogonal orthogonal rational functions.
A.~Zhedanov \cite{Zhe} pointed out the connection of
the latter continued fractions  and the generalized eigenvalue problem~\cite{Wilk}. For
applications of the $R_{\rm I}$ and $R_{\rm II}$   fractions to integrable systems, see~\cite{Spi:Zhe}. There is a more abstract approach to biorthogonality presented in Brezinski's monograph~\cite{Bre} with many applications to numerical analysis and computational mathematics.

\section{Notation}\label{section2}

Let $V_n$ be the $(2n+1)$-dimensional real vector space spanned by the
Laurent polynomials $\{z^{-n},z^{1-n},\dots, z^{n-1}, z^n\}$. Most of this
section is devoted to notation for specif\/ic spanning sets in~$V_n$.
These sets have orthogonality relations given in later sections, we
def\/ine the appropriate measures here.

We assume throughout the paper that
$x = \cos \theta$ and $z = e^{i\theta}$,
and ${\bf t}$ stands for the vector  $(t_1, t_2, t_3, t_4)$. We shall always assume that
$0<q<1$ and each $t_i$ is real with $|t_i|<1$.

The weight function for the Askey--Wilson polynomials  is given by
\begin{equation*}
w_{\text{aw}}(x;{\bf t}\,|\,q)
= \frac{\big(q,e^{2i\theta}, e^{-2i\t};q\big)_\infty}{ \prod\limits_{j=1}^4\big(t_je^{i\theta}, t_j e^{-i\t};q\big)_\infty}.
\end{equation*}
The Askey--Wilson orthogonality  \cite[(15.2.4)]{Ism05} is
\begin{gather}
\frac{1}{2\pi} \int_{-1}^1 p_m(x;\mb{t}\,|\, q)
p_n(x;\mb{t}\,|\, q)  w_{\text{aw}}(x;\mb{t}\,|\, q)  \frac{dx}{\sqrt{1-x^2}}
\nonumber\\
\qquad{} = \prod\limits_{1 \le j<k\le 4}\left(t_j t_k;q\right)_n
 \frac{\big(q,  t_1t_2 t_3 t_4q^{n-1};q\big)_n}
{(t_1t_2t_3t_4;q)_{2n}}
\mu(\mb{t}\,|\,q)  \delta_{m,n},\label{eqaworthrel}
\end{gather}
while the Askey--Wilson integral evaluation is \cite{Ask2:Wil,Ism05}
\begin{gather*}
\frac{1}{2\pi} \int_0^\pi \frac{\big(q,e^{2i\theta}, e^{-2i\t};q\big)_\infty}
{\prod\limits_{j=1}^4\big(t_je^{i\theta}, t_j e^{-i\t};q\big)_\infty}  d\theta
= \frac{(t_1t_2t_3t_4;q)_\infty}
{ \prod\limits_{1 \le j<k\le 4}  (t_j t_k ;q)_\infty}
=\mu(\mb{t}\,|\,q).
\end{gather*}

The orthogonality relations for Laurent polynomials in $V_n$ uses a small variation of the
Askey--Wilson weight, which is given by
\begin{gather*}
w_{\text{cher}}(z;{\bf t}\,|\,q) = w_{\text{aw}}(x;{\bf t}\,|\,q) \frac{(1-t_1/z)(1-t_2/z)}{1-z^{-2}}
 = \frac{\big(q,z^2, qz^{-2};q\big)_\infty}{ \prod\limits_{j=1}^4(t_j z;q)_\infty
\prod\limits_{k=3}^4(t_k/z, qt_{k-2}/z;q)_\infty}.
\end{gather*}
Even though this weight is not positive on the unit circle, and does not lead to
a positive def\/inite quadratic form, we will def\/ine a bilinear form
using $w_{\text{cher}}(z;{\bf t}\,|\,q)$,
and f\/ind orthogonal bases for~$V_n$ in Sections~\ref{section3} and~\ref{section4}.

The orthogonality relations will involve Laurent polynomials in $z$
which may be def\/ined in terms of  ${}_{4}\phi_3$ functions.  We shall need four of these functions:
$R_n$, $S_n$, $T_n$, and~$U_n$.

The f\/irst polynomial $R_n(z, {\bf t}\,|\,q)$ is an unscaled version of the Askey--Wilson polynomials
\begin{equation}
\label{eqDefRn}
R_n(z; {\bf t}\,|\,q)= {}_{4}\phi_3\left(\left. \begin{matrix}
q^{-n}, t_1t_2t_3t_4q^{n-1}, t_1z, t_1/z \\
t_1t_2,  t_1t_3,  t_1t_4
\end{matrix}  \right|q,q\right).
\end{equation}
Put $S_0(z)=0$, and for $n\ge 1$,
\begin{gather}
\label{eqDefSn}
 S_n(z; {\bf t}\,|\,q)= z(1-t_3/z)(1-t_4/z)\,
 {}_{4}\phi_3\left(\left. \begin{matrix}
q^{1-n}, t_1t_2t_3t_4q^{n}, t_1z, qt_1/z \\
qt_1t_2,  qt_1t_3,   qt_1t_4
\end{matrix}  \right|q,q\right),
\\
 \label{eqDefTn}
T_n(z; {\bf t}\,|\,q)= {}_{4}\phi_3\left(\left. \begin{matrix}
q^{-n}, t_1t_2t_3t_4q^{n-1}, t_1z, qt_1/z \\
qt_1t_2,  t_1t_3,  t_1t_4
\end{matrix}  \right|q,q\right),
\end{gather}
and $U_0(z)=0$, and for $n\ge 1$,
\begin{equation}
\label{eqDefUn}
 U_n(z; {\bf t}\,|\,q)=\frac{1}{z} (1-t_1z)(1-t_2z) R_{n-1}(z;qt_1,qt_2,t_3,t_4\,|\,q).
\end{equation}
The above representation \eqref{eqDefUn}
expresses $U_n$ as a multiple of a scaled~$R_n$. The~$S_n$ and~$T_n$
in~\eqref{eqDefSn} and \eqref{eqDefTn} can also written in terms of scaled $R_n$'s as follows
\begin{gather*}
S_n(z; {\bf t}\,|\,q)  =  z(1-t_3/z)(1-t_4/z) T_{n-1}(z; t_1, t_2, qt_3, qt_4),
\\
  T_n(z; {\bf t}\,|\,q)  =  R_n\big(q^{-1/2}z; q^{1/2} t_1, q^{1/2} t_2, q^{-1/2} t_3,
  q^{-1/2} t_4\big).
\end{gather*}
In other words
\begin{gather*}
S_n(z; {\bf t}\,|\,q)
= z(1-t_3/z)(1-t_4/z)
R_{n-1}\big(q^{-1/2}z; q^{1/2}t_1, q^{1/2}t_2, q^{1/2}t_3, q^{1/2}t_4\big).
\end{gather*}

Clearly
\begin{gather*}
\{R_0,\ldots, R_n, S_1,\ldots, S_n, T_0,\ldots, T_n, U_1,\ldots, U_n\} \subset V_n,
\end{gather*}
but this set could not form a basis for $V_n$ for $n\ge 0$. There are very specif\/ic
dependencies amongst these Laurent polynomials which follow from contiguous relations.
It follows from  Proposition~\ref{linearrelation} that for $n\ge 1$ the
subspace $W_n=\operatorname{span}\{R_n, S_n, T_n, U_n\}$ of $V_n$
is 2-dimensional. Any two of these four Laurent polynomials are
linearly independent.

We need some contiguous relations for balanced $_4\phi_3$'s. Let
\begin{equation*}
\phi(a,b,c,d;e,f,g) =  {}_{4}\phi_3\left(\left. \begin{matrix}
a, b, c, d \\
e, f, g
\end{matrix}  \right|q,q\right)
\end{equation*}
denote a ${}_4\phi_3$ which is terminating and balanced.
By $\phi_+$ we mean change $a, b, \dots, g$ by $aq, bq$,  $\dots, gq$, respectively, so
$\phi_+$  is no longer balanced. When a parameter
$\alpha$ is changed to~$\alpha q^{\pm 1}$ we write~$\phi(\alpha\pm)$.
Wilson~\cite{Wil} proved contiguous relations which imply
\begin{gather}
\label{eqcontig1}
\phi(a+, e+)-\phi = \frac{q(a-e)(1-b)(1-c)(1-d)}{(1-e)(1-eq)(1-f)(1-g)}\phi_+(e+),
\\
 \frac{a(1-f/a)(1-g/a)}{(1-f)(1-g)}\phi_+(a-) + \frac{(b-e)(1-c/e)}{(1-b)(1-c)}
\phi(d+,e+)  = \frac{(1-e)(1-bc/e)}{(1-b)(1-c)}\phi.
\label{eqcontig2}
\end{gather}

\begin{Proposition}
\label{linearrelation}
The following connection
relations hold for $n\ge 1$,
\begin{gather*}
  T_n(z;\mb{t}\,|\,q)-R_n(z;\mb{t}\,|\,q)=  \frac{qt_1\big(1-q^{-n}\big)\big(1-t_1t_2t_3t_4q^{n-1}\big)}
{(1-t_1t_2)(1-qt_1t_2)(1-t_1t_3)(1-t_1t_4)} U_n(z;\mb{t}\,|\,q),    \\
   (1-t_1t_3)(1-t_1t_4)R_n(z;\mb{t}\,|\,q) -t_1 S_n(z; {\bf t}\,|\,q)
= \frac{t_1\big(1-t_1t_2q^n\big)\big(1-t_3t_4q^{n-1}\big)}
{q^{n-1}(1-t_1t_2)(1-qt_1t_2)} U_n(z;\mb{t}\,|\,q), \\
  \frac{t_1}{(1-t_1t_3)(1-t_1t_4)}S_n(z; {\bf t}\,|\,q) - \frac{\big(1- q^n t_1t_2\big)\big(1- q^{n-1}t_3t_4\big)}
{\big(1-q^n\big)\big(1-t_1t_2t_3t_4q^{n-1}\big)} T_n(z; {\bf t}\,|\,q)  \\
 \qquad {} +   \frac{q^{n}(1-t_1t_2)(1- t_3t_4/q)}{\big(1-q^n\big)\big(1-t_1t_2t_3t_4q^{n-1}\big)}
R_n(z; {\bf t}\,|\,q)=0.
\end{gather*}
\end{Proposition}

\begin{proof} The f\/irst statement follows from \eqref{eqcontig1} with
$a=t_1/z$, $b=q^{-n}$, $c=t_1t_2t_3t_4q^{n-1}$, $d=t_1z$, $e=t_1t_2$, $f=t_1t_3$,  and
$g=t_1t_4$. The third statement follows from from~\eqref{eqcontig2} with
$a=t_1z$, $b=q^{-n}$,
$c=t_1t_2t_3t_4q^{n-1}$, $d=t_1/z$, $e=t_1t_2$, $f=t_1t_3$,  and
$g=t_1t_4$. The second statement follows from the other two.
\end{proof}

We will also use the Sears transformation \cite[(III.15)]{Gas:Rah},
\begin{gather}
  {}_{4}\phi_3\left(\left. \begin{matrix}
q^{ -n},A, B, C  \\
D, E,   F
\end{matrix} \right|q,q\right)    = \frac{(E/A, F/A;q)_n}{(E, F;q)_n}
A^n  {}_{4}\phi_3\left(\left. \begin{matrix}
q^{ -n},A, D/B, D/C\\
D,  Aq^{1-n}/E,    Aq^{1-n}/F
\end{matrix}  \right|q,q\right),\label{eqSears}
\end{gather}
where $DEF = ABC q^{1-n}$. The Sears transformation~\eqref{eqSears}
applied to the~$R_n$ and~$T_n$ functions give
\begin{gather}
\label{eqRnint2}
R_n(z; t_1,t_2,t_3, t_4\,|\,q)= \frac{t_1^n(t_2t_3, t_3t_4;q)_n}{t_3^n(t_1t_2, t_1t_4;q)_n}
R_n(z; t_3, t_2,t_1, t_4\,|\,q)
\end{gather}
and
\begin{gather}
\label{eqTnint2}
T_n(z; t_1, t_2,t_3, t_4\,|\,q)= \frac{t_1^n(t_2t_3, t_2t_4;q)_n}{t_2^n(t_1t_3, t_1t_4;q)_n}
T_n(z; t_2, t_1, t_3, t_4\,|\,q).
\end{gather}
One may also easily see from the def\/inition of $R_n$ and $S_n$ that
\begin{gather}
\label{invert}
R_n(1/z; 1/{\bf t}\,|\,1/q)=  R_n(z; {\bf t}\,|\,q),\qquad
S_n(1/z; 1/{\bf t}\,|\,1/q)=   S_n(z; {\bf t}\,|\,q)/t_3t_4.
\end{gather}

\section{Orthogonality relations}\label{section3}

The main result of this section is Theorem~\ref{orthbasis1}.

For Laurent polynomials $f(z)$ and $g(z)$ def\/ine a bilinear form by
\begin{gather*}
\langle f,g\rangle_{\text{\rm cher}}  =  \frac{1}{2\pi i} \oint_{|z|=1} f(z)
g(z)  w_{\text{\rm cher}}(z;\mb{t}\,|\, q)  \frac{dz}{z}.
\end{gather*}

\begin{Theorem}
\label{orthbasis1}
The set $\{R_0,\dots, R_n, U_1, \dots, U_n\}$
is an orthogonal basis for  $V_n$ with respect to the
symmetric bilinear form $\langle\,,\,\rangle_{\rm cher}$. Moreover, $\langle R_m,R_m\rangle_{\rm cher}$ and
 $\langle U_m,U_m\rangle_{\rm cher}$ are non-zero.
\end{Theorem}

The proof of Theorem~\ref{orthbasis1} will be accomplished in three stages.
We get  $R_n-R_m$ orthogonality for ``free''. We establish in
Propositions~\ref{eqSnorthrel} and~\ref{RTorth} orthogonalities between
$R_n$ and $T_m$, which changes into $R_n-U_m$ and
$U_n-U_m$ orthogonality.

First we establish $R_n-R_m$ orthogonality from Askey--Wilson orthogonality.
It is a restatement of~\eqref{eqaworthrel}, which was given in \cite[Lemma~6.4] {NS}.
\begin{Proposition}\label{eqSnorthrel}
 The orthogonality relation for the $R_n$ functions is
\begin{gather*}
\langle R_m,R_n\rangle_{{\text{\rm cher}}}
=\frac{(1-t_1t_2) t_1^{2n}\big(q,t_2t_3,t_2t_4,t_3t_4, t_1t_2 t_3 t_4q^{n-1};q\big)_n}
{(t_1t_2,t_1t_3,t_1t_4;q)_n(t_1t_2t_3t_4;q)_{2n}} \mu(\mb{t}\,|\,q)\delta_{m,n}.
\end{gather*}
\end{Proposition}

\begin{proof}
Rewrite the left-hand side of Proposition~\ref{eqSnorthrel} as
\begin{gather*}
\frac{1}{2\pi} \int_{-\pi}^\pi R_m\big(e^{i\t};\mb{t}\,|\, q\big)
R_n\big(e^{i\t};\mb{t}\,|\, q\big)  w_{\text{aw}}(\cos \t;\mb{t}\,|\, q)
  \frac{e^{i\t} - (t_1+t_2) + t_1t_2 e^{-i\t}}{2\sin \t}    d\t  \\
\qquad{} = \frac{1-t_1t_2}{2\pi} \int_{0}^\pi R_m\big(e^{i\t};\mb{t}\,|\, q\big)
R_n\big(e^{i\t};\mb{t}\,|\, q\big)  w_{\text{aw}}(\cos \t;\mb{t}\,|\, q)   d\t,
\end{gather*}
since $R_n\big(e^{i\t};\mb{t}\,|\, q\big)$ and $w_{\text{aw}}$ are even functions of~$\t$.
\end{proof}

To establish the $U_n-R_m$ orthogonality we f\/irst give an integral evaluation
for ``gluing'' a~specif\/ic Laurent polynomial onto the measure.
\begin{Lemma}
\label{intlemma1}
For non-negative integers $k$ and $j$ we have
\begin{gather*}
\frac{1}{2\pi i} \oint_{|z|=1}  (t_1z, qt_1/z;q)_j(t_3z, t_3/z;q)_k
  w_{\text{\rm cher}}(z;\mb{t}\,|\, q)  \frac{dz}{z}
\\
\qquad{}= \frac{\big(t_1t_2t_3t_4q^{j+k};q\big)_\infty}{\big(t_1t_2q^{j+1}, t_1t_3q^{j+k},
 t_1t_4q^j, t_2t_3q^k, t_2t_4, t_3t_4q^k;q\big)_\infty}.
\end{gather*}
\end{Lemma}
\begin{proof}
The integral is
\begin{gather*}
\frac{1}{2\pi i} \oint_{|z|=1}  \frac{\big(q,z^2,  z^{-2};q\big)_\infty}
{\big(t_1q^jz, t_2z, t_3q^kz, t_4z, q^{j}t_1/z, t_2/z, q^kt_3/z, t_4/z;q\big)_\infty}
  \frac{(1-t_2/z)\big(1-q^jt_1/z\big)}{\big(1-z^{-2}\big)}    \frac{dz}{z}.
\end{gather*}
This time the rational function in $z=e^{i\theta}$ in the integrand is
\begin{gather}
\label{ratfunc}
\frac{e^{i\theta}-\big(q^jt_1+t_2\big)+t_1t_2q^je^{-i\theta}}{2i\sin(\theta)}.
\end{gather}
The remaining part of the integrand is symmetric in $\theta$, so only
the odd part of the numerator in~\eqref{ratfunc} contributes. The Askey--Wilson
integral evaluation gives the result
\begin{gather*}
{}=\frac{1}{2}  \big(1-t_1t_2q^j\big)
\frac{1}{2\pi } \int_{-\pi}^\pi  w_{\text{aw}}\big(\cos \theta; q^jt_1, t_2, q^kt_3, t_4\,|\,q\big)   d\theta  \\
\qquad{} = \frac{\big(t_1t_2t_3t_4q^{j+k};q\big)_\infty}{\big( t_1t_2q^{j+1}, t_1t_3q^{j+k},
t_1t_4q^j, t_2t_3q^k, t_2t_4, t_3t_4q^k;q\big)_\infty}.\tag*{\qed}
\end{gather*}
\renewcommand{\qed}{}
\end{proof}

 \begin{Proposition}
 \label{RTorth}
 We have the orthogonality relation
 \begin{gather*}
\langle R_m,T_n\rangle _{{\text{\rm cher}}}
= \langle R_m,R_n\rangle_{{\text{\rm cher}}} \delta_{m,n}.
\end{gather*}
\end{Proposition}

\begin{proof} The proof is done in three stages: $m<n$, $m=n$, and $m>n$.

Let $I_{j,k}$ be the integral in Lemma~\ref{intlemma1}.  Using the explicit def\/inition of
$T_n(z;\mb{t}\,|\, q)$ we have
\begin{gather*}
\frac{1}{2\pi i} \oint_{|z|=1}  (t_3z, t_3/z;q)_k
T_n(z;\mb{t}\,|\, q)  w_{\text{cher}}(z;\mb{t}\,|\, q)  \frac{dz}{z}
= \sum_{j=0}^n \frac{\big(q^{-n}, t_1t_2t_3t_4q^{n-1};q\big)_j}
{(q, q t_1t_2, t_1t_3, t_1t_4;q)_j}q^j I_{j,k} \\
\qquad{} = \frac{(t_1t_2t_3t_4q^k;q)_\infty}
{\big(q, qt_1t_2, q^kt_1t_3, t_1t_4, q^kt_2t_3, t_2t_4, q^kt_3t_4;q\big)_\infty}
\,  {}_{3}\phi_2\left(\left. \begin{matrix}
q^{-n}, t_1t_2t_3t_4q^{n-1}, q^k t_1t_3 \\
t_1t_3,  t_1t_2t_3 t_4q^{k}
\end{matrix}  \right|q,q\right).
\end{gather*}
From the $q$-Pfaf\/f--Saalsch\"{u}tz theorem  \cite[(II.12)]{Gas:Rah}, the
sum of the above ${}_3\phi_2$ contains a factor $(q^{-k};q)_n$ which
vanishes for $k < n$. Thus if
$m<n$, then $\langle R_m,T_n\rangle_{{\text{cher}}}= 0$.

When $k=n$ the the $q$-Pfaf\/f--Saalsch\"{u}tz theorem \cite[(II.12)]{Gas:Rah}
evaluates the above integral to
\begin{gather*}
\frac{1}{2\pi i} \oint_{|z|=1}  (t_3z, t_3/z;q)_n
T_n(z;\mb{t}\,|\, q)\, w_{\text{cher}}(z;\mb{t}\,|\, q)  \frac{dz}{z}\\
\qquad{}= \frac{\big(t_1t_2t_3t_4q^{2n};q\big)_\infty  (-t_1t_3)^n
q^{\binom{n}{2}}}{\big(q^{n+1},q t_1t_2, t_1t_3, t_1t_4, t_2t_3q^n,  t_2t_4q^n,  t_3t_4q^n;q\big)_\infty}.
\end{gather*}
Using \eqref{eqRnint2}  an explicit computation shows that
Proposition~\ref{RTorth} holds for $m = n$.

Next we consider the case $m >n$. In this case consider the integral
\begin{gather*}
I=\frac{1}{2\pi i} \oint_{|z|=1}  (t_1z, q t_1/z;q)_j
R_m(z;\mb{t}\,|\, q)  w_{\text{cher}}(z;\mb{t}\,|\, q)  \frac{dz}{z}. 
\end{gather*}
After using \eqref{eqRnint2} for $R_m$, $I$ is a constant multiple of
\begin{gather*}
\sum_{k=0}^m \frac{\big(q^{-m}, t_1t_2t_3t_4q^{m-1};q\big)_k  q^k}{(q, t_1t_3, t_2t_3, t_3t_4;q)_k}
I_{j,k}.
\end{gather*}
Applying  Lemma \ref{intlemma1} we see that the integral $I$
is a constant multiple of
\begin{gather*}
{}_{3}\phi_2\left(\left. \begin{matrix}
q^{-m}, t_1t_2t_3t_4q^{m-1}, q^j t_1t_3 \\
t_1t_3,  t_1t_2t_3 t_4q^{j}
\end{matrix}  \right|q,q\right),
\end{gather*}
which again by the $q$-Pfaf\/f--Saalsch\"{u}tz theorem \cite[(II.12)]{Gas:Rah}
produces a factor $(q^{-j};q)_m$ and therefore
vanishes for $j < m$.
\end{proof}

Finally we come to the $T_n-T_m$ orthogonality, and f\/inding
the $L^2$-norms.
\begin{Proposition}\label{TTorth}
 We have the orthogonality relation
 \begin{gather*}
\langle T_m,T_n\rangle_{{\text{\rm cher}}}
= \frac{(1-t_1t_2)(1-t_3t_4/q)}{(1-t_1t_2q^n)(1-t_3t_4q^{n-1})} q^n
\langle R_m,R_n\rangle_{{\text{\rm cher}}} \delta_{m,n}.
\end{gather*}
 \end{Proposition}

\begin{proof}
We will use \eqref{eqDefTn} and \eqref{eqTnint2} to expand the $T$'s,
so consider the integral
 \begin{gather*}
 \frac{1}{2\pi i} \oint_{|z|=1}
 (t_1z, qt_1/z;q)_j   (t_2z, qt_2/z;q)_k  w_{\text{cher}}(z;\mb{t}\,|\, q)  \frac{dz}{z}\\
 \qquad{} =  \frac{1}{2\pi i} \oint_{|z|=1} (t_1z, t_1/z;q)_j   (t_2z, t_2/z;q)_k
 \frac{\big(z-q^kt_2\big)\big(1-q^jt_1/z\big)}{z-1/z}w_{\text{aw}}(z;\mb{t}\,|\,q)  \frac{dz}{z}\\
\qquad{} = \frac{(1-t_1t_2q^{k+j})}{2\pi}  \int_0^\pi w_{\text{aw}}\big(\cos \theta; t_1q^j, t_2q^k,
 t_3,t_4\,|\,q\big)d \theta \\
\qquad{}  = \frac{\big(t_1t_2t_3t_4q^{k+j};q\big)_\infty}
 {\big(q, t_1t_2q^{j+k+1}, t_1t_3q^j, t_1t_4q^j, t_2t_3q^k, t_2t_4 q^k, t_3t_4 ;q\big)_\infty}.
 \end{gather*}
Using \eqref{eqDefTn} to sum on $j$ we f\/ind
 \begin{gather*}
 \frac{1}{2\pi i} \oint_{|z|=1}
 T_n(z;\mb{t}\,|\, q)  (t_2z, qt_2/z;q)_k  w_{\text{cher}}(z;\mb{t}\,|\, q)  \frac{dz}{z}\\
 \qquad = \frac{\big(t_1t_2t_3t_4q^k;q\big)_\infty}{\big(t_1t_2q^{k+1}, t_1t_3, t_1t_4, t_2t_3q^k, t_2t_4q^k,
 t_3t_4;q\big)_\infty} \, {}_{3}\phi_2\left(\left. \begin{matrix}
q^{-n}, t_1t_2t_3t_4q^{n-1}, t_1t_2q^{k+1} \\
qt_1t_2,  t_1t_2t_3 t_4q^{k}
\end{matrix}  \right|q,q\right).
 \end{gather*}
 Again the $q$-Pfaf\/f--Saalsch\"{u}tz theorem \cite[(II.12)]{Gas:Rah}
 evaluates the ${}_3\phi_2$ as
 \begin{gather*}
 \frac{\big(q^{-k}, q^{2-n}/t_3t_4;q\big)_n}{\big(qt_1t_2, q^{1-k-n}/t_1t_2t_3t_4;q\big)_n},
\end{gather*}
which vanishes when $k < n$.  This establishes orthogonality for $m<n$, and
the $k=n$ case establishes Proposition~\ref{TTorth} for $m=n$ using  \eqref{eqTnint2}.
 \end{proof}

To complete the proof of Theorem~\ref{orthbasis1},
we switch the orthogonality results for $R_n$ and $T_n$ to results
for $R_n$ and $U_n$ using Proposition~\ref{linearrelation},
\begin{equation*}
c_nU_n= T_n-R_n, \qquad c_n=\frac{qt_1\big(1-q^{-n}\big)\big(1-t_1t_2t_3t_4q^{n-1}\big)}
{(1-t_1t_2)(1-qt_1t_2)(1-t_1t_3)(1-t_1t_4)}.
\end{equation*}
The only remaining orthogonality relation we must check is
\begin{equation*}
\langle U_n,R_n\rangle _{\rm cher}= c_n^{-1}\langle T_n-R_n,R_n\rangle_{\rm cher}
=c_n^{-1}(\langle T_n,R_n\rangle_{\rm cher}-\langle R_n,R_n\rangle_{\rm cher})=0.
\end{equation*}
The value of $\langle U_n,U_n\rangle_{\rm cher}$, which is non-zero,  is found using Proposition~\ref{TTorth},
\begin{gather*}
\langle U_n,U_n\rangle_{\rm cher}=  c_n^{-2}\langle T_n-R_n,T_n-R_n\rangle_{\rm cher}
= c_n^{-2}(\langle T_n,T_n\rangle _{\rm cher}-\langle R_n,R_n\rangle_{\rm cher})\\
\hphantom{\langle U_n,U_n\rangle_{\rm cher}}{}
=  c_n^{-2} \left(\frac{(1-t_1t_2)(1-t_3t_4/q)}{\big(1-t_1t_2q^n\big)\big(1-t_3t_4q^{n-1}\big)}
q^n-1\right)\langle R_n,R_n\rangle\\
\hphantom{\langle U_n,U_n\rangle_{\rm cher}}{}
=  c_n^{-2} \frac{\big(1-q^n\big)\big(1-t_1t_2t_3t_4q^{n-1}\big)}{\big(1-t_1t_2q^n\big)\big(1-t_3t_4q^{n-1}\big)}
\langle R_n,R_n\rangle_{\rm cher}.
\end{gather*}

\section[More orthogonal bases for $V_n$]{More orthogonal bases for $\boldsymbol{V_n}$}\label{section4}

In this section we f\/ind additional orthogonal bases for $V_n$, see Theorems~\ref{orthbasis2},~\ref{orthbasis3} and~\ref{orthbasis4}. We also
explicitly give three and four term recurrence relations for some of these bases, which are
the analogues of the recurrence relation for orthogonal polynomials.

The f\/irst result of this section is
\begin{Theorem}
\label{orthbasis2}
The set $\{T_0,\dots, T_n, S_1, \dots, S_n\}$
is an orthogonal basis for  $V_n$ with respect to the
symmetric bilinear from $\langle\,,\,\rangle_{\rm cher}$. Moreover, $\langle T_m,T_m\rangle_{\rm cher}$ and
 $\langle S_m,S_m\rangle_{\rm cher}$ are non-zero.
\end{Theorem}

\begin{proof}
Proposition~\ref{TTorth} gives orthogonality of the $T_n$'s.
By Proposition~\ref{linearrelation}
there are linear dependencies amongst $\{R_n,U_n,T_n\}$, and
$\{R_n,U_n,S_n\}$. So $\langle S_m,S_m\rangle $, $\langle T_m,T_m\rangle $, and $\langle S_m,T_m\rangle $
may all be changed to the appropriate linear combinations of~$R_n$'s and~$T_m$'s.
The results are given in the next proposition.
\end{proof}

\begin{Proposition}
\label{SSorth}
For all $n\ge1$ and  $m\ge 0$ we have the orthogonality relations
\begin{gather*}
\langle S_n,T_m\rangle_{\rm cher}   =   0, \\
t_1\langle S_n,R_m\rangle_{\rm cher}   =   (1-t_1t_3)(1-t_1t_4) \langle R_n,R_m\rangle_{\rm cher} .
\end{gather*}
For all $n\ge1$ and  $m\ge 1$ we have the orthogonality relations
\begin{gather*}
\langle S_n,S_m\rangle_{\rm cher}  =\frac{q^n(1-t_1t_3)^2(1-t_1t_4)^2(1-t_1t_2)(1-t_3t_4/q)}
{\big(q^n-1\big)\big(1-t_1t_2t_3t_4q^{n-1}\big)t_1^2}
\langle R_n,R_m\rangle_{\rm cher}.
\end{gather*}
\end{Proposition}

\begin{proof}  This follows by direct calculation using Propositions~\ref{linearrelation},~\ref{RTorth} and~\ref{TTorth}.
\end{proof}

\begin{Corollary}
\label{twoorth}
If $n\ge 1$, then for any $f\in V_{n-1}$ we have
\begin{gather*}
\langle R_n,f\rangle_{\rm cher}  = \langle S_n,f\rangle_{\rm cher} =  \langle T_n,f\rangle_{\rm cher} =  \langle U_n,f\rangle _{\rm cher} = 0.
\end{gather*}
\end{Corollary}

Using Corollary \ref{twoorth} one may explicitly give many orthogonal
bases using the Askey--Wilson polynomials~$R_n$.  For example, for
the basis
\begin{gather*}
\{R_0,\dots, R_n, T_1,\dots, T_n\}
\end{gather*}
of $V_n$,
the only pairwise non-orthogonalities involve $T_m$ and $R_m$. One
can replace these two polynomials by orthogonal linear combinations
$R_m+c_mT_m$ and $R_m+d_mT_m$. The values of the constants~$c_m$ and~$d_m$ yielding orthogonality can be chosen using
Propositions~\ref{RTorth} and~\ref{TTorth}.  The
non-symmetric Askey--Wilson polynomial of Section~\ref{section6} are one
way of accomplishing this.

Another such choice corresponds to applying the Gram--Schmidt
process to the ordered
basis $\{1,z^{-1}, z^1, \dots, z^{-n}, z^n\}$ on $V_n$ with $\langle\,,\,\rangle_{\rm cher}$.

\begin{Theorem}
\label{orthbasis3}
An orthogonal basis for $V_n$ is
$\{X_{0}, X_{-1}, X_1,\dots, X_{-n}, X_n\}$, where
\begin{gather*}
X_{-n}=  (1-t_1t_2)(1-t_1t_3)(1-t_1t_4)R_n-t_1\big(1-t_1t_2t_3t_4q^{n-1}\big)S_n, \qquad n>0,\\
X_{n}=  (1-t_1t_2)(1-t_1t_3)(1-t_1t_4)R_n+t_1^2t_2\big(1-q^{n}\big)S_n, \qquad n\ge 0.
\end{gather*}
\end{Theorem}
\begin{proof}  The choice of the coef\/f\/icients for $X_{-n}$ allows
$X_{-n}\in V_{n-1}\oplus sp\{z^{-n}\}$, while the coef\/f\/icients for
 $X_{n}$ force $\langle X_n, X_{-n}\rangle_{\rm cher}=0$ via Proposition~\ref{SSorth}.
 An explicit computation shows that all of the $L^2$-norms are non-zero.
 \end{proof}

For the ordered basis $\{1,z^{1}, z^{-1}, \dots, z^{n}, z^{-n}\}$ on $V_n$
we have an analogous result.

\begin{Theorem}
\label{orthbasis4}
An orthogonal basis for $V_n$ is
$\{Y_{0}, Y_{-1}, Y_1,\dots, Y_{-n}, Y_n\}$, where
\begin{gather*}
Y_{-n}=  q^{n-1}t_3t_4(1-t_1t_2)(1-t_1t_3)(1-t_1t_4)
R_n-t_1\big(1-t_1t_2t_3t_4q^{n-1}\big)S_n, \qquad n>0,\\
Y_{n}=  q^n(1-t_1t_2)(1-t_1t_3)(1-t_1t_4)R_n+t_1\big(1-q^{n}\big)S_n, \qquad n\ge 0.
\end{gather*}
\end{Theorem}

Orthogonal bases from Theorems \ref{orthbasis3} and~\ref{orthbasis4}
have a three term recurrence relation.

\begin{Proposition}
\label{3term}
For $n\ge 1$, there exist constants $a_n$, $b_n$ and $c_n$ such that
\begin{gather*}
zX_{-n}(z)=a_n X_{-n}(z)+b_n X_n(z) +c_n X_{n-1}(z).
\end{gather*}
\end{Proposition}

\begin{proof} Because $X_{-n}$ has no $z^n$ term, $zX_{-n}(z)\in V_n$ and
has an expansion
\begin{gather*}
zX_{-n}(z)=c_0X_0(z)=\sum_{i=1}^n \bigl(c_i X_i(z)+c_{-i}X_{-i}(z) \bigr).
\end{gather*}
By orthogonality, we can f\/ind the coef\/f\/icients
\begin{equation*}
c_k= \frac{\langle zX_{-n},X_k\rangle_{\rm cher}}{\langle X_k,X_k\rangle_{\rm cher}}
=\frac{\langle X_{-n},zX_k\rangle_{\rm cher}}{\langle X_k,X_k\rangle _{\rm cher}},
\end{equation*}
as long as $\langle X_k,X_k\rangle _{\rm cher}\neq 0$
(which is true here). The only terms to survive the orthogonality are $X_n$, $X_{-n}$, and
 $X_{n-1}$. The term $X_{1-n}$ does not survive because $zX_{1-n}\in V_{n-1}$:
\begin{gather*}
\langle zX_{-n},X_{1-n}\rangle_{\rm cher}=\langle X_{-n},zX_{1-n}\rangle_{\rm cher}=0.\tag*{\qed}
\end{gather*}
\renewcommand{\qed}{}
\end{proof}

We record two propositions for these recurrences. The coef\/f\/icients may be
found by consi\-de\-ring specif\/ic powers of~$z$ on each side, we do not give
the details of the computation.

\begin{Proposition} The constants in Proposition~{\rm \ref{3term}} are
\begin{gather*}
a_n= \frac{\big({-}t_1 - t_2 + t_1t_2(t_3 + t_4)q^{n - 1}\big)\big(1 - t_1t_2t_3t_4q^{2 n - 1}\big)}
    {\big({-}1 - t_1t_2 +  t_1t_2t_3t_4q^{n - 1} + t_1t_2q^n\big)\big(1 - t_1t_2t_3t_4q^{2 n - 2}\big)},\\
b_n=  -\frac{\big({-}t_1 - t_2 + t_1t_2(t_3 + t_4)q^{n - 1}\big)\big(1 -  t_1t_2t_3t_4q^{n - 1}\big)\big(1 - t_3t_4q^{n - 1}\big)}
{\big({-}1 - t_1t_2 + t_1t_2t_3t_4q^{n - 1} + t_1t_2q^n\big)\big(1 - t_1t_2t_3t_4q^{2 n - 2}\big)},\\
c_n=  \frac{t_1\big(1 - t_3t_4q^{n - 1}\big)\big(1 - t_2t_4q^{n - 1}\big)\big(1 - t_2t_3q^{n - 1}\big)}
   {\big(1 - t_1t_2q^{n - 1}\big)\big(1 - t_1t_2t_3t_4q^{2 n - 2}\big)}.
\end{gather*}
\end{Proposition}

\begin{Proposition}\label{2nd3term}
For $n\ge 1$ we have
\begin{gather*}
\frac{1}{z} X_n= a_n X_{-n-1}+b_n X_{n}+c_n X_{-n}, \\
a_n=  \frac{\big(1 - t_1t_2q^{n }\big)\big(1 - t_1t_3q^{n}\big)(1 - t_1t_4q^{n}\big)}
{t_1\big(1 - t_3t_4q^{n }\big)\big(1 - t_1t_2t_3t_4q^{2 n}\big)},\\
 b_n=  \frac{\big({-}t_1-t_2+t_1t_2(t_3+t_4)q^n\big)\big(1 - t_1t_2t_3t_4q^{2 n-1}\big)}
 {\big({-}1-t_1t_2+t_1t_2t_3t_4q^{n-1}+t_1t_2q^n\big)\big(1 - t_1t_2t_3t_4q^{2 n}\big)},
 \\
c_n=  - \frac{\big(1-q^n\big)\big(1 - t_1t_2q^{n }\big)\big({-}t_1-t_2+t_1t_2(t_3+t_4)q^n\big)}
{\big({-}1-t_1t_2+t_1t_2t_3t_4q^{n-1}+t_1t_2q^n\big)\big(1 - t_1t_2t_3t_4q^{2 n}\big)}.
\end{gather*}
\end{Proposition}

There are four term relations for $zX_n$ and $\frac{1}{z}X_{-n}$.

\begin{Proposition}\label{4term}
For $n\ge 1$, there exist constants $a_n$, $b_n$, $c_n$ and $d_n$ such that
\begin{gather*}
zX_{n}(z)=a_n Y_{-n-1}(z)+b_n X_n(z) +c_n Y_{-n}(z)+d_n X_{n-1}(z).
\end{gather*}
\end{Proposition}

\begin{proof} By orthogonality we have
\begin{gather*}
zX_{n}(z) =a'_n X_{-n-1}(z)+b'_n X_{n+1}(z) + c'_n X_{-n}(z)+d'_n X_{n}(z)  +e'_n X_{n-1}(z).
\end{gather*}
The sixth term $X_{1-n}(z)$ does not appear because $zX_{1-n}\in V_{n-1}$:
\begin{gather*}
\langle zX_n,X_{1-n}\rangle_{\rm cher}=\langle X_n,zX_{1-n}\rangle_{\rm cher}=0.
\end{gather*}

Since $\operatorname{span}\{X_{-n-1},X_{n+1}\}= \operatorname{span}\{Y_{-n-1},Y_{n+1}\}$, the f\/irst two terms
may be written as linear combination of $Y_{-n-1}$ and $Y_{n+1}$, But $zX_{-n}(z)$
has no $z^{-n-1}$ term, thus only $Y_{-n-1}$ appears. Because
$\operatorname{span}\{X_{-n},X_{n}\}= \operatorname{span}\{Y_{-n},X_{n}\}$ we can replace~$X_n$ by~$Y_{-n}$.
\end{proof}

\begin{Proposition} The constants in Proposition~{\rm \ref{4term}} are
\begin{gather*}
a_n= \frac{\big(1-t_1t_2q^n\big)\big(1-t_1t_3q^n\big)\big(1-t_1t_4q^n\big)
\big({-}1-t_1t_2 + t_1t_2t_3 t_4q^{n - 1}+t_1t_2q^n\big)}
    {t_1\big(1-t_3t_4q^n\big)\big(1 - t_1t_2t_3t_4q^{2 n - 1}\big)\big(1 - t_1t_2t_3t_4q^{2 n}\big)},\\
b_n=  \frac{q^n\big({-}1-t_1t_2 + t_1t_2t_3t_4q^{n - 1}+t_1t_2q^n\big)\big({-}t_3-t_4+t_2t_3t_4q^n+t_1t_3t_4q^n\big)}
{\big(1 - t_1t_2t_3t_4q^{2 n - 1}\big)\big(1 - t_1t_2t_3t_4q^{2 n}\big)},\\
c_n= - \frac{\big(1-q^n\big)\big(1 - t_1t_2q^{n}\big)\big({-}t_1-t_2+t_1t_2(t_3+t_4)q^{n-1}\big)}
   {\big(1 - t_1t_2t_3t_4q^{2 n - 2}\big)\big(1 - t_1t_2t_3t_4q^{2 n - 1}\big)},\\
d_n=   \frac{t_1\big(1-q^n\big)\big(1-t_1t_2q^n\big)\big(1-t_2t_3q^{n-1}\big)\big(1-t_2t_4q^{n-1}\big)\big(1-t_3t_4q^{n-1}\big)}
{\big(1-t_1t_2q^{n-1}\big)\big(1 - t_1t_2t_3t_4q^{2 n - 2}\big)\big(1 - t_1t_2t_3t_4q^{2 n-1}\big)}.
\end{gather*}
\end{Proposition}

\begin{Proposition}
\label{2nd4term}
For $n\ge 1$ we have
\begin{gather*}
\frac{1}{z} X_{-n}= a_n X_{-n-1}+b_n X_{-n}+c_n Y_{n}+d_n Y_{n-1}, \\
a_n=  \frac{\big(1 - t_1t_2q^{n }\big)\big(1 - t_1t_3q^{n}\big)\big(1 - t_1t_4q^{n}\big)
\big(1-t_3t_4q^{n-1}\big)\big(1-t_1t_2t_3t_4q^{n-1}\big)}
{t_1\big(1 - t_3t_4q^{n }\big)\big(1 - t_1t_2t_3t_4q^{2 n-1}\big)\big(1 - t_1t_2t_3t_4q^{2 n}\big)},\\
b_n=  \frac{q^{n-1}\big({-}t_3-t_4+(t_1+t_2)t_3t_4q^{n-1}\big)
\big({-}1-t_1t_2+t_1t_2t_3t_4q^{n-1}+t_1t_2q^n\big)}
 {\big(1 - t_1t_2t_3t_4q^{2 n-2}\big)\big(1 - t_1t_2t_3t_4q^{2 n-1}\big)},
 \\
c_n=  - \frac{\big(1 - t_3t_4q^{n-1 }\big)\big(1 - t_1t_2t_3t_4q^{n-1}\big)\big({-}t_1-t_2+t_1t_2(t_3+t_4)q^n\big)}
{\big(1 - t_1t_2t_3t_4q^{2 n-1}\big)\big(1 - t_1t_2t_3t_4q^{2 n}\big)},\\
d_n=  \frac{t_1\big(1 - t_2t_3q^{n-1 }\big)\big(1 - t_2t_4q^{n-1}\big)\big(1 - t_3t_4q^{n-1}\big)
\big({-}1-t_1t_2+t_1t_2t_3t_4q^{n-1}+t_1t_2q^n\big)}
{\big(1 - t_1t_2q^{n-1 }\big)\big(1 - t_1t_2t_3t_4q^{2 n-2}\big)\big(1 - t_1t_2t_3t_4q^{2 n-1}\big)}.
\end{gather*}
\end{Proposition}

The two bases of Theorems~\ref{orthbasis3} and~\ref{orthbasis4} can be related,
using~\eqref{invert}. To explicitly show the parameter dependence we write
$X_{n}(z;{\mb{t}}\,|\,q)$ for $X_n$ and $Y_{n}(z;{\mb{t}}\,|\,q)$ for~$Y_n$.

\begin{Proposition}\label{invertXY}  If $n\ge 1$, we have
\begin{gather*}
X_{-n}(1/z;1/{\mb{t}}\,|\,1/q)=   -Y_{-n}(z;{\mb{t}}\,|\,q)/\big(t_1^3t_2t_3^2t_4^2q^{n-1}\big),\\
X_{n}(1/z;1/{\mb{t}}\,|\,1/q)=   -Y_{n}(z;{\mb{t}}\,|\,q)/\big(t_1^3t_2t_3t_4q^{n}\big).
\end{gather*}
\end{Proposition}

Thus three term relations can be explicitly given for $zY_n$ and $\frac{1}{z}Y_{-n}$, and
four term relations for $zY_{-n}$ and $\frac{1}{z}Y_{n}$.

\section[The operators $A_0$ and $A_1$]{The operators $\boldsymbol{A_0}$ and $\boldsymbol{A_1}$}\label{section5}

Motivated by operators previously given in Noumi--Stokman \cite{NS}, in this section
we def\/ine linear transformations $A_0$ and $A_1$ on
$V_n$. We explicitly f\/ind their actions on the possible bases $R_m$,
$S_m$, $T_m$, and $U_m$. We use the self-adjointness of
these operators to give an alternative proof of Theorems~\ref{orthbasis1},~\ref{orthbasis2},~\ref{orthbasis3}, and~\ref{orthbasis4}. We describe the eigenvalue problem whose solutions
are the non-symmetric Askey--Wilson polynomials.
We explicitly give the creation operators of
Noumi--Stokman~\cite{NS}
via the three and four term recurrence relations of Section~\ref{section4}.
These results are in~\cite{NS}, but we include them for completeness.

For a Laurent polynomial $f\in V_n$, def\/ine
\begin{gather*}
(A_0f)(z) = \frac{(1-t_3/z)(1-t_4/z)}{(1-q/z^2)}(f(q/z)-f(z)),
\end{gather*}
and
\begin{gather*}
(A_1f)(z) = \frac{(1-t_1z)(1-t_2z)}{(1-z^2)}(f(1/z)-f(z)).
\end{gather*}

By a direct computation it is easy to see that each~$A_i$ is a~linear transformation from~$V_n$ to~$V_n$.

\begin{Proposition} \label{prop5.1}
Each $A_i:V_n\rightarrow V_n$ is self-adjoint, that is for any
$f,g\in V_n$, we have
\begin{gather*}
\langle A_0f,g\rangle _{\rm cher}=\langle f,A_0g\rangle _{\rm cher}, \qquad
\langle A_1f,g\rangle _{\rm cher}=\langle f,A_1g\rangle_{\rm cher}.
\end{gather*}
\end{Proposition}

\begin{proof} It is clear that $\langle A_0f,g\rangle_{\rm cher}$ equals
\begin{gather*}
\frac{1}{2\pi i} \oint_{|z|=1} \frac{\big(q, z^2, 1/z^{2};q\big)_\infty}{ \prod\limits_{j=1}^4(t_jz, t_j /z;q)_\infty}
\frac{(1-t_1/z)(1-t_2/z)}{1-z^{-2}}  \\
\quad \qquad {}\times  \frac{(1-t_3/z)(1-t_4/z)}{\big(1-q/z^2\big)}\left[f(q/z)g(z)-f(z)g(z)\right]
\frac{dz}{z}   \\
\qquad{}=    \frac{1}{2\pi i} \oint_{|z|=1} \frac{\big(q, z^2, q^2/z^{2};q\big)_\infty}{ \prod\limits_{j=1}^4(t_jz, qt_j /z;q)_\infty}
   \left[f(q/z)g(z)-f(z)g(z)\right]
\frac{dz}{z}
\end{gather*}
This a dif\/ference of two integrals. In the f\/irst we let $z \to q/z$.
A calculation shows that  the f\/irst integral becomes
\begin{gather*}
 \frac{1}{2\pi i} \oint_{|z|=q} \frac{\big(q, z^2, q^2/z^{2};q\big)_\infty}
 { \prod\limits_{j=1}^4(t_jz, qt_j /z;q)_\infty}
 f(z)g(q/z)
\frac{dz}{z}.
\end{gather*}
We now deform the above contour to the circular contour $|z| =1$
because the integrand has no poles between the contours since
each $|t_i|<1$. The proof of $\langle A_1f,g\rangle_{\rm cher}=\langle f,A_1g\rangle_{\rm cher}$ is similar and will be omitted.
\end{proof}

Note that in the proof of Proposition~\ref{prop5.1} we only
used the fact that~$f$ and~$g$ are analytic in $q \le |z| \le 1/q$.

\begin{Proposition}
\label{operaction}
The action of the operators~$A_0$ and~$A_1$ on the
functions $R_n$, $T_n$, $S_n$, and~$U_n$ is given by
\begin{gather*}
  A_0(R_n)= \frac{t_1\big(q^{-n}-1\big)\big(1-t_1t_2t_3t_4q^{n-1}\big)}
{(1-t_1t_2)(1-t_1t_3)(1-t_1t_4)} S_n=\alpha_n S_n,  \qquad
 A_1(R_n)= 0,\\
  A_0(T_n)=0,  \qquad A_1(T_n) = \frac{qt_1\big(q^{-n}-1\big)\big(1-t_1t_2t_3t_4q^{n-1}\big)}
{(1-qt_1t_2)(1-t_1t_3)(1-t_1t_4)} U_n=\beta_n U_n, \\
  A_0(S_n)=(-1+t_3t_4/q)S_n,  \qquad  A_1(S_n)=
\frac{q^{1-n}\big(1-t_1t_2q^{n}\big)\big(1-t_3t_4q^{n-1}\big)}
{(1-qt_1t_2)}U_n,\\
 A_0(U_n)= q^{-1}(1-qt_1t_2)S_n, \qquad  A_1(U_n)=(-1+t_1t_2) U_n,
\end{gather*}
respectively.
\end{Proposition}
\begin{proof} These follow by direct term by term application of $A_0$
and $A_1$ to the def\/ining series.
\end{proof}

These  relations give another proof of the orthogonality relations in Theorems~\ref{orthbasis1}
 and~\ref{orthbasis2}, for example
\begin{gather*}
0=\langle A_0(T_n), R_m\rangle_{\rm cher}= \langle T_n, A_0(R_m)\rangle_{\rm cher}=
\alpha_m \langle T_n, S_m\rangle_{\rm cher},
\\
0=\langle A_1(R_n), U_m\rangle_{\rm cher}= \langle R_n, A_1(U_m)\rangle_{\rm cher}=
\beta_m \langle R_n, U_m\rangle_{\rm cher}.
\end{gather*}

From Proposition~\ref{operaction} the operators $A_0$ and $A_1$ act on the
2-dimensional space $W_n$ spanned by $\{R_n,S_n,U_n, T_n\}$. An elementary
computation yields the next two theorems.

\begin{Theorem}
\label{firsteigen}
If $n\ge 1$, the eigenvalues of $(A_1-t_1t_2I)\big(A_0-q^{-1}t_3t_4I\big)$ on $W_n$
are~$q^{-n}$ and $t_1t_2t_3t_4q^{n-1}$, with corresponding eigenvectors
$X_{-n}(z;{\bf{t}}\,|\,q)$ and $Y_n(z;{\bf{t}}\,|\,q).$
\end{Theorem}

\begin{Theorem}
If $n\ge 1$, the eigenvalues of
$\big(A_0-q^{-1}t_3t_4I\big)(A_1-t_1t_2I)$ on $W_n$
are $q^{-n}$ and $t_1t_2t_3t_4q^{n-1}$, with corresponding eigenvectors
$Y_{-n}(z;{\bf{t}}\,|\,q)$ and $X_n(z;{\bf{t}}\,|\,q)$.
\end{Theorem}

These two theorems give alternative proofs of
Theorems~\ref{orthbasis3} and~\ref{orthbasis4}. For example
\begin{gather*}
q^{-n}\langle X_{-n}, X_n\rangle_{\rm cher} =  \langle (A_1-t_1t_2I)\big(A_0-q^{-1}t_3t_4I\big)X_{-n}, X_n\rangle_{\rm cher}\\
\hphantom{q^{-n}\langle X_{-n}, X_n\rangle_{\rm cher}}{}
=  \langle X_{-n}, \big(A_0-q^{-1}t_3t_4I\big)(A_1-t_1t_2I)X_n\rangle_{\rm cher}  \\
\hphantom{q^{-n}\langle X_{-n}, X_n\rangle_{\rm cher}}{}
=  t_1t_2t_3t_4q^{n-1}\langle X_{-n}, X_n\rangle_{\rm cher}
\end{gather*}
implying that $0=\langle X_{-n}, X_n\rangle_{\rm cher}.$

Our next goal is to def\/ine ladder maps $\mathbb{S}_0$ and $\mathbb{S}_1$ such that
$\mathbb{S}_1$ maps $Y_n$ to $X_{-n-1}$, and
$\mathbb{S}_0$ maps $X_{-n}$ to $Y_{n}$.
 If such maps are found, then~$Y_n$ is an $n^{th}$
iterate of their composition~$\mathbb{S}_0\mathbb{S}_1$
applied to $Y_0$, see
Theorems~\ref{iterate1} and~\ref{iterate2}.

Maps on $W_n$ that interchange the bases elements $\{X_{-n},Y_n\}$ and
$\{X_{n},Y_{-n}\}$ can be def\/ined using the commutators of~$A_i$ with
\begin{gather*}
Y=(A_1-t_1t_2I)\big(A_0-q^{-1}t_3t_4I\big) \qquad \text{and} \qquad
{\mathbb{Y}}=\big(A_0-q^{-1}t_3t_4I\big)(A_1-t_1t_2I).
\end{gather*}

From Proposition~\ref{operaction} the actions of $Y$, ${\mathbb{Y}}$,
$A_0$ and $A_1$ are
explicitly known on the 2-dimensional space $W_n=\operatorname{span}\{R_n,S_n,U_n, T_n\}$.
An explicit computation yields the next two propositions.

\begin{Proposition}\label{commutator1}
For $n\ge 1$ we have
\begin{gather*}
  [Y,A_1](X_{-n})= t_1t_2q^{-n}\big(1-t_3t_4q^{n-1}\big)\big(1-t_1t_2t_3t_4q^{n-1}\big) Y_n,\\
[{\mathbb{Y}},A_0](Y_{-n})= -q^{-1}t_3t_4\big(1-t_3t_4q^{n-1}\big)\big(1-t_1t_2t_3t_4q^{n-1}\big) X_n,
 \end{gather*}
and if $n\ge 0$,
\begin{gather*}
[Y,A_1](Y_{n})= q^{-n}\big(1-q^{n}\big)\big(1-t_1t_2q^{n}\big) X_{-n},\\
[{\mathbb{Y}},A_0](X_{n})= -q^{-2n}\big(1-q^{n}\big)\big(1-t_1t_2q^{n}\big) Y_{-n}.
 \end{gather*}
\end{Proposition}

\begin{Proposition} For $n\ge 1$ we have
\begin{gather*}
  [Y,A_1](R_{n})=  \frac{t_1q^{1 - 2 n}\big(1 - q^n\big)\big(1 - t_3t_4q^{n - 1}\big)
\big(1 - t_1t_2q^{n}\big)\big(1 - t_1t_2t_3t_4q^{n - 1}\big)}
{(1 - t_1t_2)(1 - qt_1t_2)/(1 - t_1t_3)/(1 - t_1t_4)} U_n,\\
[Y,A_1](U_{n})=  q^{-1}t_2(1 - t_1t_4)(1 - t_1t_2)(1 - qt_1t_2) (1 - t_1t_3) R_n,\\
 [Y,A_1](S_{n})=  -q^{-n}t_2\big(1-q^nt_1t_2\big)(1 - t_1t_3)(1 - t_1t_4)\big(1 - t_3t_4q^{n-1}\big) R_n\\
\hphantom{[Y,A_1](S_{n})=}{}
+q^{1-2n}\frac{\big(1-q^n\big)\big(1-t_1t_2q^n\big)\big(1-t_3t_4q^{n-1}\big)\big(1- t_1t_2t_3t_4q^{n - 1}\big)}
{(1-t_1t_2)(1-qt_1t_2)} U_n,\\
[Y,A_1](T_{n})=  -q^{-n}t_1t_2\big(1 - t_1t_2t_3t_4q^{n-1}\big) R_n\\
\hphantom{[Y,A_1](T_{n})=}{}
 +q^{1-2n}t_1\frac{\big(1-q^n\big)\big(1-t_1t_2q^n\big)\big(1-t_3t_4q^{n-1}\big)\big(1- t_1t_2t_3t_4q^{n - 1}\big)}
{(1-t_1t_2)(1-qt_1t_2)(1-t_1t_3)(1-t_1t_4)} U_n.
\end{gather*}
\end{Proposition}

We see from Proposition \ref{commutator1}
that  the commutator $\mathbb{S}_0=[Y,A_1]$ does map
$X_{-n}$ to~$Y_{n}$. To def\/ine an appropriate~$\mathbb{S}_1$, we need
another pair of operators.

We def\/ine operators $\hat T_1, \hat T_2:V_n\rightarrow V_{n+1}$ by
\begin{gather*}
(\hat T_1f)(z)=\frac{(A_1 f)(z)+f(z)}{z}, \qquad
(\hat T_2f)(z)=z((A_0 f)(z)+f(z)).
\end{gather*}

\begin{Proposition}
\label{hatT1action}
The action of  $\hat T_1$ on $Y_n$ and $X_{-n}$ is given by
\begin{gather*}
\hat T_1(Y_n)=  \frac{\big(1 - q^n t_1 t_2\big) \big(1 - q^nt_1 t_3\big) \big(1 - q^n t_1 t_4\big)}
{t_1 \big(1 - q^n t_3 t_4\big) \big(1 - q^{2 n} t_1 t_2 t_3 t_4\big)}X_{-n-1}
 + \frac{(t_1 + t_2 - q^nt_1t_2(t_3 + t_4))}{(1 - q^{2 n}t_1t_2t_3t_4)}Y_{n},
\end{gather*}
if $n \ge 0$  and
\begin{gather*}
\hat T_1(X_{-n})=  q^{n - 1} \frac{t_1 t_2 \big({-}t_4 + t_3 \big(-1 + q^{n - 1} (t_1 + t_2) t_4\big)\big)}
{\big({-}1 + q^{2 (n - 1)} t_1 t_2 t_3 t_4\big) }X_{-n}\\
\hphantom{\hat T_1(X_{-n})=}{}
+ \frac{ t_1^2 t_2 \big(1 - q^{n - 1} t_2 t_3\big) \big(1 - q^{n - 1} t_2 t_4\big)
\big(1 - q^{n - 1}t_3 t_4\big)}
{\big(1 - q^{n - 1} t_1 t_2\big) \big({-}1 + q^{2 (n - 1)} t_1 t_2 t_3 t_4\big)} Y_{n-1},
\end{gather*}
if $n\ge 1$, respectively.
\end{Proposition}

\begin{proof}[Proof (sketch).] Note that $\hat T_1= \frac{1}{z}(A_1+I)$. The action of $A_1$ on any
element of $W_n$ is given by Proposition~\ref{operaction}.
The action of division by $z$ is given by
Propositions~\ref{2nd3term},~\ref{2nd4term}, and their $Y$ versions. Thus
$\hat T_1(Y_n)$ and $\hat T_1(X_{-n})$ may be given explicitly.
\end{proof}

From Theorem~\ref{firsteigen} we know that $Y_{n}$ and $X_{-n}$ are
eigenfunctions of $Y$,
and can obtain Proposition~\ref{T1comm} from Proposition~\ref{hatT1action}.

\begin{Proposition}\label{T1comm}
For $n\ge 0$ we have
\begin{gather*}
  [Y,\hat T_1](Y_{n})= -\frac{\big(1-t_1t_2q^{n}\big)\big(1-t_1t_3q^{n}\big)\big(1-t_1t_4q^{n}\big)}
{t_1q^{n+1}\big(1-t_3t_4q^{n}\big)} X_{-n-1},\\
[Y,\hat T_1](X_{-n-1})= -t_1^2t_2\frac{\big(1-t_2t_3q^{n}\big)\big(1-t_2t_4q^{n}\big)\big(1-t_3t_4q^{n}\big)}
{q^{n+1}\big(1-t_1t_2q^{n}\big)} Y_n.
\end{gather*}
\end{Proposition}

Finally we def\/ine $\mathbb{S}_1=[Y,\hat T_1]$, which maps $Y_n$ to $X_{-n-1}$.
The next theorem is Proposition~9.3 in~\cite{NS}.
Recall that $Y_0=(1-t_1t_2)(1-t_1t_3)(1-t_1t_4)$ is a constant.

\begin{Theorem}
\label{iterate1}
We have
\begin{gather*}
c_n Y_n=
(\mathbb{S}_0 \mathbb{S}_1)^n(Y_0), \qquad
d_n X_{-n-1}= \mathbb{S}_1 (\mathbb{S}_0\mathbb{S}_1)^{n}(Y_0).
\end{gather*}
where
\begin{gather*}
c_n=  (-t_2)^n q^{-n(n+1)} (t_1t_2,t_1t_3,t_1t_4,t_1t_2t_3t_4;q)_n,\\
d_n=  -(-t_2)^n q^{-(n+1)^2} (t_1t_2,t_1t_3,t_1t_4;q)_{n+1}(t_1t_2t_3t_4;q)_n/t_1\big(1-t_3t_4q^n\big).
\end{gather*}
\end{Theorem}

Completely analogous results hold for
\begin{gather*}
\hat T_2, \qquad \mathbb{S}_2=[{\mathbb{Y}},A_0],  \qquad\mathbb{S}_3=[{\mathbb{Y}},\hat T_2],
\end{gather*}
which we state below. Note that Proposition \ref{commutator1} shows that $\mathbb{S}_2$ interchanges
$Y_{-n}$ and $X_n$.

\begin{Proposition} For $n\ge 0$ we have
\begin{gather*}
\hat T_2(X_n)=  \frac{t_3+t_4-t_3t_4q^n(t_1+t_2)}
{\big(1 - q^{2 n} t_1 t_2 t_3 t_4\big)}X_{n}  + \frac{\big(1 - q^{n} t_1 t_2\big) \big(1 - q^{n } t_1 t_3\big)
\big(1 - q^{n}t_1 t_4\big)}{t_1q^n\big(1-q^nt_3t_4\big)\big(1 - q^{2 n}t_1t_2t_3t_4\big)}Y_{-n-1},\\
\hat T_2(Y_{-n-1})=  t_1t_3t_4q^{n } \frac{\big(1 - q^{n} t_2 t_3\big) \big(1 - q^{n } t_2 t_4\big)
\big(1 - q^{n}t_3 t_4\big)}
{\big(1 - q^{n}t_1 t_2\big)\big(1 - q^{2 n} t_1 t_2 t_3 t_4\big) }X_{-n}\\
\hphantom{\hat T_2(Y_{-n-1})=}{}
+ q^nt_3t_4\frac{\big(t_2+t_1-t_1t_2q^n(t_3+t_4)\big)}
{ \big(1 - q^{2 n } t_1 t_2 t_3 t_4\big)} Y_{-n-1}.
\end{gather*}
\end{Proposition}

The map $\mathbb{S}_3=[{\mathbb{Y}},\hat T_2]$ interchanges $X_n$ and $Y_{-n-1}$.

\begin{Proposition} For $n\ge 0$ we have
\begin{gather*}
  [{\mathbb{Y}},\hat T_2](X_{n})=  \frac{\big(1-t_1t_2q^{n}\big)\big(1-t_1t_3q^{n}\big)\big(1-t_1t_4q^{n}\big)}
{t_1q^{2n+1}\big(1-t_3t_4q^{n}\big)} Y_{-n-1},\\
[{\mathbb{Y}},\hat T_2](Y_{-n-1})= q^{-1}t_1t_3t_4\frac{\big(1-t_2t_3q^{n}\big)\big(1-t_2t_4q^{n}\big)(1-t_3t_4q^{n}\big)}
{\big(1-t_1t_2q^{n}\big)} X_n.
\end{gather*}
\end{Proposition}

Recall that $X_0=(1-t_1t_2)(1-t_1t_3)(1-t_1t_4)$ is a constant.

\begin{Theorem}
\label{iterate2}
We have
\begin{gather*}
e_n X_n= (\mathbb{S}_2 \mathbb{S}_3)^n(X_0),
\qquad f_n Y_{-n-1}= \mathbb{S}_3 (\mathbb{S}_2\mathbb{S}_3)^{n}(X_0),
\end{gather*}
where
\begin{gather*}
e_n=  (-t_3t_4/t_1)^n q^{-n(n+1)} (t_1t_2,t_1t_3,t_1t_4,t_1t_2t_3t_4;q)_n,\\
f_n=  (-t_3t_4/t_1)^n q^{-n(n+1)} (t_1t_2,t_1t_3,t_1t_4;q)_{n+1}(t_1t_2t_3t_4;q)_n/t_1q^{2n+1}\big(1-t_3t_4q^n\big).
\end{gather*}
\end{Theorem}

Since any two elements from $\{X_n,X_{-n},Y_n,Y_{-n}\}$ form a basis for $W_n$,
any element of $W_n$ has several such formulas.

\section[Non-symmetric Askey-Wilson polynomials]{Non-symmetric Askey--Wilson polynomials}\label{section6}

The non-symmetric Askey--Wilson polynomials may be def\/ined as the
eigenfunctions given in Theorem \ref{firsteigen}, namely
\begin{equation}
\{Y_0,Y_1,\ldots, Y_n, \dots\}\cup \{X_{-1},X_{-2},\ldots, X_{-n}, \dots \}.
\nonumber
\end{equation}
This set of Laurent polynomials also satisfy an biorthogonality relation which
follows from Theorems~\ref{orthbasis3} and~\ref{orthbasis4}.

Def\/ine a new bilinear form on $V_n$ by
\begin{gather*}
\langle f,g\rangle_{\text{cher}'}  =  \frac{1}{2\pi i} \oint_{|z|=1}f(z) g(1/z)
w_{\text{cher}}(z;\mb{t}\,|\, q)  \frac{dz}{z}.
\end{gather*}
Put
\begin{gather*}
P_m(z,{\bf{t}}\,|\, q)=
\begin{cases}
Y_m(z,{\bf{t}}\,|\,q) \  & \text{if} \ m\ge 0,\\
X_m(z,{\bf{t}}\,|\,q) \ & \text{if} \ m<0,
\end{cases}
\end{gather*}
and $P_m'(z,{\bf{t}}\,|\,q))=P_m(z,1/{\bf{t}}\,|\,1/q).$

The biorthogonality relation with respect to this bilinear form \cite{M}, \cite[\S~6.7]{NS} is
\begin{Corollary}
\label{nonsymmorth}
For $m\neq n$, $\langle P_m,P_n'\rangle _{{\rm cher}'}=0$.
\end{Corollary}

\begin{proof} We use Proposition~\ref{invertXY} to change the bilinear form
$\langle\,,\,\rangle_{{\rm cher}'}$ to $\langle\,,\,\rangle_{\rm cher}$. For example, if $m, n< 0$, $m\neq n$, then
\begin{gather*}
\langle X_m, X_n'\rangle_{{\rm cher}'}=-\langle X_m, Y_n \rangle_{\rm cher}/\big(t_1^3t_2t_3^2t_4^2q^{n-1}\big)=0,
\end{gather*}
while if $m<0$ and $n\ge 0$
\begin{gather*}
\langle X_m, Y_n'\rangle_{{\rm cher}'}=-\langle X_m, X_n\rangle _{\rm cher}t_1^3t_2t_3t_4q^{n}=0,
\end{gather*}
The only term which requires checking is for $n>0$,
\begin{gather*}
\langle Y_n, X_{-n}'\rangle_{{\rm cher}'}=- \langle Y_n, Y_{-n}\rangle_{\rm cher}/\big(t_1^3t_2t_3^2t_4^2q^{n-1}\big)=0
\end{gather*}
and
\begin{gather*}
\langle X_{-n}, Y_{n}'\rangle_{{\rm cher}'}=-\langle X_{-n}, X_{n}\rangle_{\rm cher}t_1^3t_2t_3^2t_4^2q^{n}=0
\end{gather*}
by Theorems~\ref{orthbasis3} and~\ref{orthbasis4}.
\end{proof}

The $L^2$-norms may also be found.

\begin{Corollary} We have for $n>0$,
\begin{gather*}
\langle Y_n, Y_n'\rangle_{{\rm cher}'}=
-q^{n} \langle R_n,R_n \rangle_{\rm cher}/\big(t_1^3t_2t_3t_4\big)\\
\hphantom{\langle Y_n, Y_n'\rangle_{{\rm cher}'}=}{}
 \times \frac{(1 - t_1 t_2) \big(1 -  t_1 t_2q^n\big) (1 - t_1 t_3)^2 (1 - t_1 t_4)^2
\big(1 -  t_1 t_2 t_3 t_4q^{2 n-1}\big)}{\big(1 -  t_1 t_2 t_3 t_4q^{n-1}\big)},
\\
\langle X_{-n}, X_{-n}'\rangle_{{\rm cher}'}= -t_1^3t_2t_3^2t_4^2q^{1-n}\langle R_n,R_n\rangle_{\rm cher}\\
\hphantom{\langle X_{-n}, X_{-n}'\rangle_{{\rm cher}'}=}{}
 \times \frac{(1 - t_1 t_2) \big(1 -  t_3 t_4q^{n-1}\big) (1 - t_1 t_3)^2 (1 - t_1 t_4)^2
\big(1 -  t_1 t_2 t_3 t_4q^{2 n-1}\big)}{\big(1 - q^{n}\big)} .
\end{gather*}
\end{Corollary}

\begin{proof} Using Proposition~\ref{invertXY} we have
\begin{gather*}
\langle Y_n, Y_n'\rangle_{{\rm cher}'}=-q^n\langle Y_n, X_n\rangle_{\rm cher}/\big(t_1^3t_2t_3t_4\big),
\end{gather*}
and $\langle Y_n, X_n\rangle_{\rm cher}$ may be found using Theorems~\ref{orthbasis3},~\ref{orthbasis4} and Proposition~\ref{SSorth}. The proof for
$\langle X_{-n}, X_{-n}'\rangle_{{\rm cher}'}$ is similar.
\end{proof}

Koornwinder and Bouzef\/four~\cite{KB} gave a positive def\/inite
orthogonality relation for the non-symmetric Askey--Wilson polynomials.

Macdonald \cite[(6.6.8), p.~167]{M} gives an expression for
$P_m(z;\mb{t}\,|\,q)$ as linear combination of $R_m(z;{\mb{t}}\,|\,q)$
and $T_m(z;{\mb{t}}\,|\,q)$. This is equivalent to the Noumi--Stokman
expression using Proposition~\ref{linearrelation}.

Another \cite[\S~6.8]{NS} biorthogonal pair of bases with
respect to $\langle\,,\,\rangle_{{\rm cher}'}$ is given by the bases
\[
\{R_0,\ldots,R_n, U_1,\ldots, U_n \} \qquad \text{and} \qquad
\{R'_0,\ldots,R_n', U_{1}',\ldots, U_{n}' \}.
\]
The $U_n(z;\mb{t}\,|\,q)$ are the anti-symmetric Askey--Wilson polynomials.
However since
\begin{gather*}
R_n(1/z;1/{\bf{t}}\,|\,1/q)= R_n(z;{\bf{t}}\,|\,q), \qquad U_n(1/z;1/{\bf{t}}\,|\,1/q)=
U_n(z;{\bf{t}}\,|\,q)/t_1t_2,
\end{gather*}
this orthogonality is equivalent to Theorem~\ref{orthbasis1}.

\section[Askey-Wilson orthogonality]{Askey--Wilson orthogonality}\label{section7}

One may ask if there are orthogonality relations for a basis of $V_n$
using the Askey--Wilson weight. By symmetry, one could take
\begin{gather*}
\{R_0,\ldots, R_n, (z-1/z),(z-1/z)R_1,\ldots, (z-1/z)R_{n-1}\}.
\end{gather*}
We include a partial result in this direction.

 \begin{Theorem}
 We have the orthogonality relation
\begin{gather*}
\frac{1}{2\pi i} \oint_{|z|=1} R_m(z;\mb{t}\,|\, q)
T_n(z;\mb{t}\,|\, q)  w_{\text{aw}}(x;\mb{t}\,|\, q) \frac{dz}{z}\\
\qquad{} =  \frac{\big({-}t_1^2\big)^n q^{\binom{n}{2}} \big(1+q^n\big)\big(t_1t_2t_3t_4q^{2n};q\big)_\infty}{2\prod\limits_{1\le j <k \le 4}(t_jt_k ; q)_\infty}
 \prod_{2\le j <k \le 4} (t_jt_k ; q)_n   \delta_{m,n},
\end{gather*}
for $m \ge n$.
 \end{Theorem}
\begin{proof}
Consider the integral
  \begin{gather*}
\frac{1}{2\pi i} \oint_{|z|=1} R_m(z;\mb{t}\,|\, q)
(t_1z, qt_1/z;q)_k  w_{\text{aw}}(x;\mb{t}\,|\, q)  \frac{dz}{z},
\end{gather*}
for $0\le k \le n$.
The integral vanishes when $k=0$ by \eqref{eqaworthrel}, so we assume $k>0$. Write
$(t_1z, qt_1/z;q)_k$ as $(1-t_1z)(1-t_1q^k/z)(qt_1z, qt_1/z;q)_{k-1}$. Therefore the above integral is
\begin{gather*}
\frac{1}{2\pi }\int_{-\pi}^\pi R_m\big(e^{i\theta};\mb{t}\,|\, q\big)
\big(qt_1e^{i\theta}, qt_1e^{-i\theta};q\big)_{k-1}
   \big[1- t_1\big(1+q^k\big) \cos \theta + t_1^2q^k\big]
w_{\text{aw}}(x;\mb{t}\,|\, q) d\theta,
\end{gather*}
which vanishes for $k < m$ by~\eqref{eqaworthrel}. If $k =m$ the above integral is
  \begin{gather*}
\frac{1+q^n}{4\pi  }\int_{-\pi}^\pi R_m\big(e^{i\theta};\mb{t}\,|\, q\big)
\big(t_1e^{i\theta}, t_1e^{-i\theta};q\big)_{n}
w_{\text{aw}}(x;\mb{t}\,|\, q) d\theta \\
\qquad{} = \frac{1+q^n}{4\pi  } \frac{(q, t_1t_2, t_1, t_3, t_1t_4;q)_n}
{q^n\big(q^{-n}, t_1t_2t_3t_4q^{n-1};q\big)_n}\int_{-\pi}^\pi R_m^2\big(e^{i\theta}\big)
 w_{\text{aw}}(x;\mb{t}\,|\, q) d\theta,
\end{gather*}
and the result follows from \eqref{eqasexpl}, \eqref{eqaworthrel} and \eqref{eqDefRn}.
\end{proof}

\section{Recurrences}\label{section8}

In this section we record three term
recurrence relations satisf\/ied by $R_n$, $S_n$, $T_n$ and $U_n$.
 The three term recurrence relation for $R_n$ is the Askey--Wilson recurrence relation
 \cite[(3.1.4)]{Koe:Swa}
  \begin{gather*}
  [z+1/z] R_n(z; {\bf t}\,|\,q)   = A_n^{(r)}R_{n+1}(z; {\bf t}\,|\,q)  +
  C_n^{(r)}R_{n-1}(z; {\bf t}\,|\,q) \\
 \hphantom{[z+1/z] R_n(z; {\bf t}\,|\,q)   =}   {} +\big[{-}A_n^{(r)}- C_n^{(r)}+ t_1+ 1/t_1\big] R_n(z; {\bf t}\,|\,q),
  \end{gather*}
 where
\begin{gather*}
 A_n^{(r)} =  \frac{\big(1-t_1t_2t_3t_4q^{n-1}\big) \prod\limits_{j=2}^4\big(1-t_1t_jq^n\big)}
 {t_1\big(1-t_1t_2t_3t_4q^{2n-1}\big)\big(1-t_1t_2t_3t_4q^{2n}\big)},
 \\
 C_n^{(r)} = \frac{t_1\big(1-q^n\big)\prod\limits_{2\le j < k \le 4}\big(1-t_j t_kq^{n-1}\big)}
 {\big(1-t_1t_2t_3t_4q^{2n-2}\big)\big(1-t_1t_2t_3t_4q^{2n-1}\big)}.
 \end{gather*}
The three term  recurrence relation for $T_n$  is
  \begin{gather*}
   [z+q/z] T_n(z; {\bf t}\,|\,q)   = A_n^{(t)} T_{n+1}(z; {\bf t}\,|\,q)  +
  C_n^{(t)} T_{n-1}(z; {\bf t}\,|\,q) \\
 \hphantom{[z+q/z] T_n(z; {\bf t}\,|\,q)   =}{}
    + \big[{-}A_n^{(t)}- C_n^{(t)} + q t_1+ 1/t_1\big] T_n(z; {\bf t}\,|\,q),
  \end{gather*}
with $n >0$,   where
 \begin{gather*}
 A_n^{(t)} =
 \frac{\big(1-t_1t_2t_3t_4q^{n-1}\big)\big(1-t_1t_2q^{n+1}\big)  \prod\limits_{j=3}^4\big(1-t_1t_jq^n\big)}
 {t_1\big(1-t_1t_2t_3t_4q^{2n-1}\big)\big(1-t_1t_2t_3t_4q^{2n}\big)},
 \\
 C_n^{(t)} = \frac{qt_1\big(1-q^n\big)\big(1-t_3t_4q^{n-2}\big)\prod\limits_{j=3}^4\big(1-t_2 t_jq^{n-1}\big)}
 {\big(1-t_1t_2t_3t_4q^{2n-2}\big)(1-t_1t_2t_3t_4q^{2n-1}\big)}.
 \end{gather*}
 On the other hand the three term recurrence relation for $U_n$ is
  \begin{gather*}
   [z+q/z] U_n(z; {\bf t}\,|\,q)   = A_n^{(u)} U_{n+1}(z; {\bf t}\,|\,q)  +
  C_n^{(u)} U_{n-1}(z; {\bf t}\,|\,q) \\
 \hphantom{[z+q/z] U_n(z; {\bf t}\,|\,q)   =}{}    + \big[{-}A_n^{(u)} - C_n^{(u)}  + q t_1+ q^{-1}/t_1\big] U_n(z; {\bf t}\,|\,q),
  \end{gather*}
with $n >1$,   where
 \begin{gather*}
 A_n^{(u)} =  \frac{\big(1-t_1t_2t_3t_4q^{n}\big)\big(1-t_1t_2q^{n+1}\big) \prod\limits_{j=3}^4\big(1-t_1t_jq^n\big)}
 {qt_1\big(1-t_1t_2t_3t_4q^{2n-1}\big)\big(1-t_1t_2t_3t_4q^{2n}\big)},
 \\
 C_n^{(u)} = \frac{qt_1\big(1-q^{n-1}\big) \big(1-t_1t_2q^n\big)\big(1-t_3t_4q^{n-2}\big)}
 {\big(1-t_1t_2t_3t_4q^{2n-2}\big)\big(1-t_1t_2t_3t_4q^{2n-1}\big)}
  \prod_{j=3}^4\big(1-t_1 t_jq^{n-1}\big)\big(1-t_2 t_jq^{n-1}\big).
 \end{gather*}
 Finally the three term recurrence relation for $S_n$ is
  \begin{gather*}
\big[q^{-1/2}z+q^{1/2}/z\big] R_n(z; {\bf t}\,|\,q)   = A_n^{(s)}R_{n+1}(z; {\bf t}\,|\,q)  +
  C_n^{(s)}R_{n-1}(z; {\bf t}\,|\,q) \\
 \hphantom{\big[q^{-1/2}z+q^{1/2}/z\big] R_n(z; {\bf t}\,|\,q)   =}{}    +\big[{-}A_n^{(s)}- C_n^{(s)}+ t_1+ 1/t_1\big] R_n(z; {\bf t}\,|\,q),
  \end{gather*}
  and the coef\/f\/icients are given by
  \begin{gather*}
 A_n^{(s)} =  \frac{\big(1-t_1t_2t_3t_4q^{n}\big) \prod\limits_{j=2}^4\big(1-t_1t_jq^n\big)}
 {t_1\big(1-t_1t_2t_3t_4q^{2n-1}\big)\big(1-t_1t_2t_3t_4q^{2n}\big)},
 \\
 C_n^{(s)} = \frac{t_1\big(1-q^n\big)\prod\limits_{2\le j < k \le 4}\big(1-t_j t_kq^{n-1}\big)}
 {\big(1-t_1t_2t_3t_4q^{2n-2}\big)\big(1-t_1t_2t_3t_4q^{2n-1}\big)}.
 \end{gather*}

 \section{Asymptotics}\label{section9}

 Consider a general balanced terminating ${}_4\phi_3$,
\begin{gather*}
 {}_{4}\phi_3\left(\left. \begin{matrix}
q^{-n}, Aq^{n-1}, B,C \\
D,  E,   F
\end{matrix}  \right|q,q\right), \qquad \text{with}\quad ABC = DEF.
\end{gather*}
Since the ${}_4\phi_3$ is symmetric in $B$ and $C$ we may assume
\[
|B| \le |C|.
\]
 Ismail and Wilson \cite{Ism:Wil} determined the large degree
 asymptotics of the Askey--Wilson polynomials. When $|B| < |C|$ their result is
 \[
 {}_{4}\phi_3\left(\left. \begin{matrix}
q^{-n}, Aq^{n-1}, B,C \\
D,  E,   F
\end{matrix}  \right|q,q\right)
=C^n  \frac{(B, D/C, E/C, F/C;q)_\infty}{(B/C, D, E, F;q)_\infty}  \big[1+\O\big(q^{n/2}\big)\big].
\]
 On the other hand when $|B|= |C|$ we let $B/C = e^{2i\theta}$ and in
 this case the Ismail--Wilson asymptotic result is
 \begin{gather*}
  {}_{4}\phi_3\left(\left. \begin{matrix}
q^{-n}, Aq^{n-1}, B,C \\
D,  E,   F
\end{matrix}  \right|q,q\right) =  \big[1+\O\big(q^{n/2}\big)\big]  \\
 \qquad {} \times  \left[\frac{C^n(B, D/C, E/C, F/C;q)_\infty}{(B/C, D, E, F;q)_\infty} +
 \frac{B^n(C, D/B, E/B, F/B;q)_\infty}{(C/B, D, E, F;q)_\infty}\right] .
 \end{gather*}

 We now record the asymptotics of $R_n$, $S_n$, $T_n$ and $U_n$.

 \begin{Theorem}
 \label{asympt}
 For $|z| = 1$ the following  asymptotic formulas hold
\begin{gather*}
  \prod_{j=2}^4 (t_1 t_j;q)_\infty R_n(z; {\bf t}) =\left(\frac{t_1}{z}\right)^n
\frac{\prod\limits_{j=1}^4(t_j z;q)_\infty}
{\big(z^2;q\big)_\infty}\big[1+\O\big(q^{n/2}\big)\big] \\  
\hphantom{\prod_{j=2}^4 (t_1 t_j;q)_\infty R_n(z; {\bf t}) =}{}
  + (t_1z)^n \frac{\prod\limits_{j=1}^4(t_j/ z;q)_\infty}
{\big(1/z^2;q\big)_\infty}\big[1+\O\big(q^{n/2}\big)\big],   \\
  \frac{\big(qz^{-2};q\big)_\infty \prod\limits_{j=2}^4(qt_1t_j;q)_\infty}
{(qt_1/z, qt_2/z, t_3/z, t_4/z;q)_\infty}S_n(z; {\bf t}) = t_1^{n-1}z^n \big[1+\O\big(q^{n/2}\big)\big], \\ 
  \frac{\big(qz^{-2};q\big)_\infty(qt_1t_2;q)_\infty \prod\limits_{j=3}^4(t_1t_j;q)_\infty}
{ \prod\limits_{j=1}^2(qt_j/z, t_{j+2}/z;q)_\infty}T_n(z; {\bf t}) = (t_1z)^n \big[1+\O\big(q^{n/2}\big)\big],  \\ 
    \frac{\prod\limits_{j=2}^4(qt_1t_j;q)_\infty}{(1-qt_1t_2)} U_n(z; {\bf t}) =
\frac{(qt_1)^{n-1}}{z^n} \frac{ \prod\limits_{j=1}^4(t_j z;q)_\infty}
{\big(z^2;q\big)_\infty}[1+\O(q^{n/2})] \\ 
\qquad{} +  (q t_1)^{n-1} z^{n-2}\frac{ (1-t_1z)(1-t_2z)\prod\limits_{j=1}^4(t_j/ z;q)_\infty}
{\big(1/z^2;q\big)_\infty (1-t_1/z)(1-t_2/z)}\big[1+\O\big(q^{n/2}\big)\big].
\end{gather*}
\end{Theorem}

\section{Discrete orthogonality}\label{section10}

The Askey--Wilson polynomials have a discrete orthogonality,
the $q$-Racah orthogonality \cite{Ask:Wil, Koe:Swa},
when $t_1t_j=q^{-N}$, for some
$j=2,3$ or $4$, and for some
positive integer $N$. In this section we give the corresponding
discrete orthogonality for the Laurent polynomials. This may be done
by deforming the contours, we do not give the details.

The $q$-Racah measure for the Askey--Wilson polynomials is
purely discrete and based upon the
very well poised terminating ${}_6\phi_5$ evaluation \cite[(II.20)]{Gas:Rah}
\begin{gather}
\sum_{k=0}^N
\frac{(t_1^2,t_1t_2,t_1t_3,t_1t_4;q)_k}
{(q,qt_1/t_2,qt_1/t_3,qt_1/t_4;q)_k}
  \frac{1-t_1^2q^{2k}}{1-t_1^2}
\left( \frac{q}{t_1t_2t_3t_4}\right)^k \nonumber\\
\qquad{} =  \frac{\big(qt_1^2,qt_1^2/t_2t_3,qt_1^2/t_2t_4,qt_1^2/t_3t_4;q\big)_\infty}
{(qt_1/t_2,qt_1/t_3,qt_1/t_4,q/t_1t_2t_3t_4;q)_\infty}.\label{vwp65}
\end{gather}

Using
\begin{gather*}
\frac{1-t_1^2q^{2k}}{1-t_1^2}=
\frac{1}{\big(1-t_1^2\big)(1-t_1t_2)}\big({-}t_1t_2\big(1-q^k\big)\big(1-q^kt_1/t_2\big)+
\big(1-t_1^2q^{k}\big)\big(1-t_1t_2q^k\big)\big)
\end{gather*}
one can rewrite the left side of the ${}_6\phi_5$ evaluation~\eqref{vwp65}
as a sum of two ${}_4\phi_3$'s
\begin{gather*}
\sum_{k=1}^N \frac{\big(qt_1^2,qt_1t_2;q\big)_{k-1}
(t_1t_3,t_1t_4;q)_k}{(qt_1/t_2,q;q)_{k-1}(qt_1/t_3,qt_1/t_4;q)_k }
\left( \frac{q}{t_1t_2t_3t_4}\right)^k (-t_1t_2)
\\
\qquad{}+\sum_{k=0}^N \frac{\big(qt_1^2,qt_1t_2,t_1t_3,t_1t_4;q\big)_k }
{(qt_1/t_2,qt_1/t_3,qt_1/t_4,q;q)_k}
\left( \frac{q}{t_1t_2t_3t_4}\right)^k.
\end{gather*}

\begin{Definition}
If $t_1t_j=q^{-N}$ for some positive
integer $N$, the Racah bilinear form on $V_N$ is def\/ined by
\begin{gather*}
\langle f,  g\rangle_{\rm Racah} = \sum_{k=1}^N \frac{\big(qt_1^2,qt_1t_2;q\big)_{k-1}
(t_1t_3,t_1t_4;q)_k}{(qt_1/t_2,q;q)_{k-1}(qt_1/t_3,qt_1/t_4;q)_k }
\left( \frac{q}{t_1t_2t_3t_4}\right)^k (-t_1t_2)f\big(t_1q^k\big)g\big(t_1q^k\big)\\
\hphantom{\langle f,  g\rangle_{\rm Racah} =}{}
+ \sum_{k=0}^N \frac{\big(qt_1^2,qt_1t_2,t_1t_3,t_1t_4;q\big)_k }
{(qt_1/t_2,qt_1/t_3,qt_1/t_4,q;q)_k}
\left( \frac{q}{t_1t_2t_3t_4}\right)^k f\big(q^{-k}/t_1\big)g\big(q^{-k}/t_1\big).
\end{gather*}
\end{Definition}

\begin{Theorem}
\label{Racah}
If $t_1t_j=q^{-N}$ for some positive integer $N$ and $j=3$ or $4$, then
\begin{gather*}
\{R_0,\dots, R_N, U_1, \dots, U_N\} , \qquad \{T_0,\dots, T_N, S_1, \dots, S_N\},\\
 \{X_{-N},\dots, X_0, X_1, \dots, X_N\} \qquad \text{and} \qquad
\{Y_{-N},\dots, Y_0, Y_1, \dots, Y_N\}
\end{gather*}
are Racah-orthogonal bases for $V_N$.
\end{Theorem}

If $t_1t_2=q^{-N},$ the Racah bilinear form contains a zero term: the $k=N$ term
in the second sum is $0$ due to $(qt_1t_2;q)_k.$ In this case $U_N$ and $T_N$
are not well-def\/ined. A slight modif\/ication of Theorem~\ref{Racah} can be formulated, we
do not give it here.

\subsection*{Acknowledgements}
The research of Mourad E.H.~Ismail is partially supported  by Research Grants Council of Hong Kong under  contracts \#101410 and \#101411 and by King Saud University, Riyadh
through grant DSFP/MATH 01.

\pdfbookmark[1]{References}{ref}
\LastPageEnding

\end{document}